\documentclass{siamart0516}

\usepackage{lipsum}
\usepackage{amsfonts}
\usepackage{graphicx}
\usepackage{epstopdf}
\usepackage{algorithmic}
\usepackage[sort,compress]{cite}
\ifpdf
  \DeclareGraphicsExtensions{.eps,.pdf,.png,.jpg}
\else
  \DeclareGraphicsExtensions{.eps}
\fi

\newcommand{\TheTitle}
{Iterative importance sampling algorithms for parameter estimation} 
\newcommand{\TheAuthors}
{Morzfeld, Day, Grout, Pau, Finsterle, Bell}

\headers{Iterative importance sampling algorithms}{\TheAuthors}

\title{{\TheTitle}\thanks{Submitted to the editors DATE.
\funding{
This material is based upon work supported by the 
U.S.~Department of Energy, Office of Science,
Office of Advanced Scientific Computing Research, 
Applied Mathematics Program under contract DE-AC02-05CH11231, 
by the National Science Foundation under grant DMS-1619630,
and by the Alfred P. Sloan Foundation through a Sloan Research Fellowship awarded to MM.
This research used resources of the National Energy Research Scientific Computing Center, a DOE Office of Science User Facility supported by the Office of Science of the U.S. Department of Energy under Contract No. DE-AC02-05CH11231.
}}}

\author{
  Matthias Morzfeld\thanks{University of Arizona, Tucson, AZ
    (\email{mmo@math.arizona.edu}).}
  \and
  Marcus S. Day\thanks{Lawrence Berkeley National Laboratory
    Berkeley, CA (\email{msday@lbl.gov},
    \email{gpau@lbl.gov},
    \email{safinsterle@lbl.gov},
    \email{jbbell@lbl.gov}).}
  \and
  Ray W. Grout\thanks{National Renewable Energy Laboratory, Golden, CO
    (\email{Ray.Grout@nrel.gov}).}
  \and
  George Shu Heng Pau\footnotemark[3]  
  \and
  Stefan A. Finsterle\footnotemark[3] 
  \and 
  John B. Bell\footnotemark[3] 
}

\usepackage{amsopn}

\ifpdf
\hypersetup{
  pdftitle={\TheTitle},
  pdfauthor={\TheAuthors}
}
\fi

\definecolor{darkred}{rgb}{0.7,0,0}
\newcommand{\Red}[1]{{\color{darkred}#1}}
\begin{document}

\maketitle

\begin{abstract}
In parameter estimation problems one computes
a posterior distribution over uncertain parameters
defined jointly by a prior distribution, a model, and noisy data.
Markov Chain Monte Carlo (MCMC) is often used for the numerical solution of such problems.
An alternative to MCMC is importance sampling,
which can exhibit near perfect scaling with the number of cores
on high performance computing systems
because samples are drawn independently.
However, finding a suitable proposal distribution is a challenging task.
Several sampling algorithms have been proposed over the past years
that take an iterative approach to constructing a proposal distribution.
We investigate the applicability of such algorithms
by applying them to two realistic and challenging test problems, 
one in subsurface flow, and one in combustion modeling.
More specifically, we implement importance sampling algorithms that 
iterate over the mean and covariance matrix of Gaussian or multivariate $t$-proposal distributions.
Our implementation leverages massively parallel computers,
and we present strategies to initialize the iterations
using ``coarse'' MCMC runs or Gaussian mixture models.
\end{abstract}



\section{Introduction}
Predicting the behavior of complex physical systems is a key requirement in science and engineering. 
One approach to predicting the behavior of such systems is through high fidelity simulation. 
However, for many systems, model uncertainties limit predictive capability. 
The idea in parameter estimation is to use experimental data 
to reduce model uncertainties and, therefore,
improve predictions of overall system behavior as follows.
One represents uncertainties in the model
by assuming that model parameters are random variables
with given ``prior'' distributions,
and one represents mismatch between model and data by a ``likelihood''.  
By Bayes' rule, the prior and likelihood define a posterior distribution
that describes how well the parameters describing the system can be determined from the available data  
-- which aspects are tightly bounded and which are less precise.
The posterior distribution is thus regarded as the solution of a parameter estimation problem 
\cite{CH13,Stuart10,Tarantola05}.

The central theme of this paper is to present a Monte-Carlo (MC) sampling method
that is efficient for parameter estimation problems for complex systems, 
and well suited to extreme-scale (parallel) computer architectures.
We first formulate the Bayesian parameter estimation problem and 
review the literature about Monte Carlo approaches to its numeral solution
(section~2).
We then motivate our use of iterative importance sampling algorithms (ISA)
and present its implementation on high performance computing systems,
as well as effective strategies to initialize the iterations (section~3).
ISA has been discussed in the literature extensively over the past decades (see, e.g.,
\cite{OhBerger89,OhBerger92,CappeEtAl04}),
but it has not yet been applied to ``realistic'' test problems with complex systems.
Our main goal in this work is to demonstrate that ISA is indeed applicable in such situations
by applying our ISA implementation to two test cases: 
estimation of permeabilities in a subsurface flow problem (section~4),
and estimation of kinetic parameters for combustion simulation (section~5).

\section{Background and problem formulation}
\label{sec:ParameterEstimationProblem}
We consider a set of parameters $\theta$
of a numerical model $\mathcal{M}(\theta)$.
Here $\theta$ is an $n_\theta$-dimensional vector,
and the model $\mathcal{M}: R^{n_\theta} \to R^{n_z}$ 
maps the parameter vector to an $n_z$ dimensional vector of outputs~$z$.
In the examples discussed below, $\theta$ represents collection of reaction rate parameters and third-body
coefficients that define a combustion model,
or $\theta$ defines a permeability field of a subsurface flow problem
(see sections~\ref{sec:Subsurface} and~\ref{sec:Combustion} below).
Uncertainty in the model is represented by uncertainty in the parameters,
i.e.,  $\theta$ is a random vector with prior probability distribution $p_\theta$.
The outputs of the model can be measured, 
and additional uncertainty in the approximation of the measurements
are represented by a random variable $v$ with given distribution $p_v$.
We assume that this uncertainty is additive,
\begin{equation}
\label{eq:ModelToData}
	z = \mathcal{M}(\theta)+v,
\end{equation}
however this assumption is not central to our approach.
The above equation defines the likelihood $p(z\vert\theta) = p_v(z-\mathcal{M}(\theta))$.
By Bayes' rule, prior and likelihood define the posterior distribution
\begin{equation*}
	p(\theta\vert z) \propto p_\theta(\theta)p(z\vert\theta),
\end{equation*}
which describes the probability of the parameters $\theta$ given the data $z$.
The posterior distribution allows for a complete description of the uncertainty
in the parameters given the data. 
In particular, one can compute least-squares optimal estimates
by computing the posterior mean
as well as ``error bars'' based on posterior covariances.

If the prior and the errors $v$ are Gaussian,
and if, in addition,
the model $\mathcal{M}$ is linear, 
then the posterior is also Gaussian.
In this case, it is sufficient to compute its mean and covariance.
This can be done, for example, by minimizing the negative logarithm
of the posterior distribution
\begin{equation*}
	F(\theta) =  -\log p_\theta(\theta) - \log p(z\vert\theta) + C.
\end{equation*}
The minimizer of $F$ is the posterior mean,
and the Hessian, $H_{ij} = \partial^2 F/\partial\theta_i\partial\theta_j$, $i,j=1,\dots,n_\theta$
is the inverse of the posterior covariance matrix \cite{Tarantola05}.
Linear algorithms have been developed that are highly efficient for large-scale problems of this type
\cite{FlathEtAl11,BuiEtAl13}.

When the model is nonlinear, 
the posterior distribution is no longer Gaussian,
even if the prior distribution for the parameters and 
the model errors, $v$, remain Gaussian.
In this case, minimizing the negative log-posterior $F(\theta)$
results in the most likely state given the data, 
which is called ``posterior mode''.
However, there may be more than one posterior mode
(local maxima of the posterior distribution),
in which case it is difficult to argue that a (local) posterior mode 
is a meaningful solution of the parameter estimation problem.

\subsection{Markov Chain Monte Carlo for inverse problems}
Markov Chain Monte Carlo (MCMC) methods
are routinely used for solving nonlinear parameter estimation problems.
The basic idea is to choose a specific proposal distribution to propose a ``move'',
i.e., a new set of parameters. 
This move is accepted or rejected based on a suitable
accept\slash reject criterion \cite{Kalos86}.
The chain then consists of a number of correlated samples
and converges in the sense that averages over the samples 
converge to expected values with respect to the posterior distribution
as the ensemble size goes to infinity.
Thus, MCMC provides a more complete description of the uncertainty
than computing modes of the posterior distribution.
In practice, the beginning of an MCMC chain is often regarded
as a ``burn-in'' period, and discarded in the subsequent analysis.
In addition, if the correlation ``time''
of the chain is too large then prohibitively long chains are needed to effectively sample the posterior.
The various MCMC algorithms in the literature
differ in their proposal distributions and approaches to reduce burn-in and auto-correlation time.

There are several ``classical'' MCMC algorithms,
e.g., Metropolis-Hastings \cite{Kalos86},
and many of them are known to be 
slow in high-dimensional and highly nonlinear problems,
see, e.g., \cite{Hanson01,GoodmanWeare10}.
Recent MCMC methods that address these problems
include affine invariant ensemble samplers~\cite{GoodmanWeare10,Hammer}, adaptive MCMC \cite{HaarioEtAl01,HaarioEtAl06},
and differential evolution MCMC~\cite{Braak2006,Braak2008}.
Another class of  MCMC algorithms relies on 
using local geometry of the posterior to construct proposals.
Algorithms of this type include the
Metropolis adjusted Langevin algorithm (MALA) 
\cite{Roberts98},
Hamiltonian MCMC,
\cite{Duane87,Neal11}
and
Riemann manifold MCMC
\cite{Girolami11}.
The stochastic Newton MCMC algorithm
\cite{MartinEtAl12,PetraEtAl14}
also falls within this class.
In stochastic Newton, one first finds the posterior mode by 
minimizing the negative log-posterior,
$\min_\theta F(\theta)$,
and then starts an MCMC chain using low-rank approximations
of Hessians, evaluated along the chain.
Adjoints are used to speed up gradient and Hessian computations.
Another related method is ``randomize-then-optimize'' (RTO) \cite{RTO},
where a stochastic version of the negative log-posterior
optimization problem is solved (repeatedly) to generate samples.

In many situations, the model $\mathcal{M}$ is a discretization of a 
partial differential equation (PDE)
and the unknown parameter is a discretization of a field.
There is a large literature addressing the fact that
MCMC methods for such problems should have a well-defined
limit as the mesh is refined, see, e.g., \cite{Stuart10,CuiEtAl16b,BJMS15,CRSW13}.
Moreover, indicators of computational requirements (burn-in or autocorrelation times)
of such ``function space'' MCMC algorithms,
should not increase as the dimension increases when the mesh is refined.
In fact these properties should be invariant
under mesh-refinement, once grid convergence is reached. 
For that reason, 
function space MCMC methods are also often called ``mesh-independent''
or ``dimension-independent''.

Computational requirements of MCMC can be reduced by
solving a related, but easier problem.
This idea is implemented by reduced-order modeling:
rather than using the model $\mathcal{M}$ for MCMC,
one replaces the model by a simpler, ``reduced-order'' version of it.
Since the reduced-order model (ROM) is simpler, 
it is computationally less expensive,
and using it therefore speeds up MCMC.
In this context, one can also reduce the dimension of the parameter vector $\theta$,
the dimension of the state of the underlying numerical model, $\mathcal{M}$,
or both, see, e.g., \cite{LieberManEtAl16,CuiEtAl16,SpantiniEtAl15,CuiEtAl14}.
Reducing the state dimension is usually achieved by
finding a simplified approximation to the mathematical model~$\mathcal{M}$.
This reduced model however may depend on the same number of parameters as the ``full'' model.
To reduce this number of parameters
one can determine a low-dimensional subspace of parameter combinations
which are constrained by data, 
as in the likelihood-informed methods \cite{CuiEtAl16,CuiEtAl14}.
In problems where the parameters are spatial fields,
e.g., permeability fields in subsurface flow problems,
one can use principal component analysis\slash Karhunen-Loeve expansions \cite{MTWC15},
to define this field in terms of a small number of dominant modes
(see also section~\ref{sec:Subsurface}, where Kriging is used to reduce the parameter dimension).
In either case, MCMC by reduced-order models requires that
the ROM does not introduce large errors, 
or else the posterior distribution defined by the ROM
may be very different from the posterior distribution defined by $\mathcal{M}$.
Such ideas were made precise in the context of polynomial chaos expansions in \cite{LuEtAl15}.
To reduce error due to model reduction one can introduce a second stage in each Metropolis step
as in \cite{CF05}, or one can continually refine the ROM during MCMC sampling
as in \cite{CMPS16}.
Nonetheless, ROM approaches 
are only effective if the model $\mathcal{M}$ indeed has a low-dimensional representation
that can be discovered by ROM techniques.

Another approach to reducing computational cost of MCMC 
involves the use optimal transport maps, see, e.g.,
\cite{Reich15,Tarek12,Dean14,Recca11}. 
The idea is to construct a transport map that converts prior samples to posterior samples.
More generally, the map should transform a reference distribution
to the posterior distribution \cite{CLMMT16}.
Construction of an optimal map
requires solving an infinite dimensional optimization problem,
which is expensive,
but approximations to an optimal map
may be more readily computable and can be used
to ``speed up'' MCMC. 
Work in this direction is underway \cite{ParnoEtAl16,ParnoMarzouk16}.

\subsection{Importance sampling for inverse problems}
An alternative to MCMC is importance sampling. 
Importance sampling, just as MCMC, requires a proposal distribution~$q$.
But rather than using the proposal distribution for defining a Markov chain,
one draws samples from it and attaches to each sample a weight
\begin{equation*}
	w(\theta) \propto \frac{p(\theta\vert z)}{q(\theta)}.
\end{equation*}
The weighted samples form an empirical estimate of the posterior distribution
in the sense that weighted averages over the samples converge to 
expected values as the ensemble size goes to infinity, see, e.g., \cite{CH13,Owen}.
Moreover, the posterior distribution and the proposal distribution
need only be known up to a multiplicative constant (as in MCMC),
because the weights can be self-normalized such that their sum is one, 
i.e., $w_i \leftarrow w_i/\sum_{j=1}^{N_e} w_j$.
Advantages of importance sampling compared to MCMC are that there is no burn-in time,
and the weighted samples are independent.
The latter means in particular that the samples and weights can be 
computed independently,
which can lead to enormous computational advantages on massively parallel computer architectures.

However, the efficiency of importance sampling hinges on constructing 
a suitable proposal distribution $q$.
The reason is that, even though all samples are independent, 
each sample carries a different weight and, therefore,
the samples are not equally important.
For an effective importance sampling method 
all samples should be similarly important,
i.e., the variance of the weights should be small.
Indeed, a heuristic ``effective sample size'' can be defined 
\cite{GLM15,Weare2013,LiuChen98,DoucetEtAl01,ArulampalamEtAl02,Bergman99}
based on the relative variance $var(w)/E(w)^2$:
\begin{equation}
\label{eq:R}
	N_\text{eff} = \frac{N}{R}, \quad R =\frac{var(w)}{E(w)^2}+1 = \frac{E(w^2)}{E(w)^2},
\end{equation}
The effective sample size describes the size of
an unweighted ensemble of size $N_\text{eff}$
that is equivalent to a weighted ensemble of size $N$.
As an illustration, suppose that one weight is close to one, 
which implies that all other self-normalized weights are close to zero.
In this case, $R\approx N$, so that the effective sample size is one.
If this happens, the sampling method is said to have ``collapsed''
\cite{BengtssonEtAl08}.
In contrast, if the proposal distribution is equal to the posterior distribution,
the variance of the weights is zero, so that $R$ is unity.
For effective sampling, one needs an $R$ that is ``not too large'',
where what large is depends mainly on how much computing power one is willing to spend.

The collapse of sequential Monte Carlo methods for time-dependent problems,
called particle filters \cite{DoucetEtAl01},
has been studied extensively in the meteorological literature,
and it was shown that,
for a certain class of problems,
$N_\text{eff}\to 1$ as the ``effective dimension'' of the problem increases
\cite{CM13,SnyderEtAl08,BickelEtAl08,BengtssonEtAl08,Snyder11,SBM15,CLMMT16}.
The effective dimension can depend on the dimension of the data vector $z$,
the dimension of the parameter vector $\theta$,
as well as the distributions of errors in the prior and likelihood
and may increase when the dimension of the parameter vector $\theta$ increases,
however, at a slower rate.
Moreover, recent work suggests that particle filters
may be efficient even if the effective dimension is large,
provided that the filtering problem has a certain ``sparse'' structure~\cite{RVH2015,MHS17}.

Perhaps due to these numerical difficulties,
applications of importance sampling are rare,
but importance sampling schemes have been used
for some time-dependent ``filtering problems''
(see, e.g., \cite{Chorin09,ChorinEtAl10,MorzfeldEtAl11,AtkinsEtAl13,Reich15,ArulampalamEtAl02,Weare2013})
and a few (simple) parameter estimation test problems 
(see, e.g., \cite{MTWC15,BJMS15}).

\subsection{Iterative importance sampling}
\begin{figure}[t]
\centering
\includegraphics[width=1\textwidth]{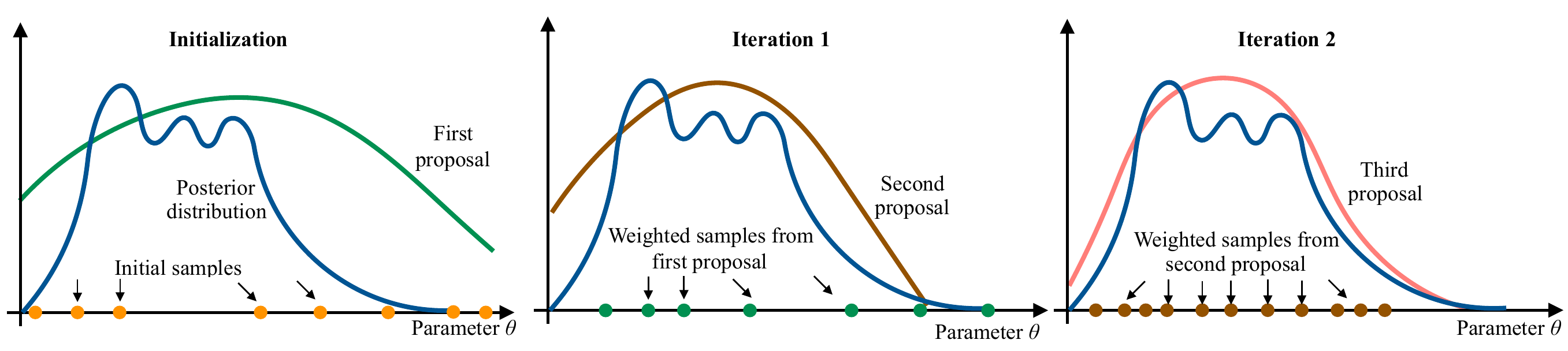}
\caption{
Illustration of iterative construction of a proposal distribution.
}
\label{fig:IllustrateAlgorithm}
\end{figure}

The basic idea in iterative importance sampling algorithms (ISA)
is to iteratively improve the proposal distribution.
The iteration starts with a given proposal distribution $q^0$,
which is used to draw $N_e^0$ weighted samples $\{\theta_j^1,w^1_j\}$, $j=1,\dots, N_e^0$.
These samples are then used to define a new proposal distribution $q^1$,
and the process is repeated.
The rationale for an iterative approach is to discover,  
via iteration, the global structure of the posterior distribution
without becoming stuck in local minima.
A successful iteration is illustrated in Figure~\ref{fig:IllustrateAlgorithm}.

ISA are discussed extensively in the literature
and are sometimes called ``adaptive importance sampling''
or ``population Monte Carlo''.
The earliest account we could find is the applications paper by Buchner 
\cite{Bucher88}, where iterative sampling is used to compute failure probabilities in structural dynamics. 
The proposal distribution at each step is Gaussian and the mean and covariance of
the Gaussian are updated at each step of the iteration. 
The paper also presents initialization strategies for the iteration that work well for the application.
A more general version of ISA is described in \cite{OhBerger92,OhBerger89}.
Several stopping rules are presented,
and consistency of the method is established,
in the sense that expected values computed by ISA are asymptotically equal to 
expected values computed with respect to the posterior distribution.
A method closely related to ISA is ``population Monte Carlo'' \cite{CappeEtAl04,Iba00},
or ``non-parametric importance sampling'' \cite{Zhang96}.
In these methods, the ensemble at each step undergoes Kernel smoothing,
i.e., the proposal distribution at iteration $k$ is
\begin{equation*}	
	q^k = \sum_{j=1}^{N_e} \alpha^k_j K_j(\theta,\theta_j^k),
\end{equation*}
where the sum of the $\alpha_j$ is one.
Here, one can adapt the weights of the mixture $\alpha_j^k$ \cite{DoucEtAl07},
or adapt the weights and the Kernels \cite{CappeEtAl08,GivensEtAl96} in each iteration.
A recent application of population Monte Carlo in geosciences is presented in \cite{Stordal15}.
The scaling of the computational requirements for one step of population Monte Carlo 
with respect to the dimension of the parameter vector
is discussed in a simplified setting in \cite{LeeEtAl11}
(see also our conclusions in section~\ref{sec:Conclusions}).
We wish to point out that the adaptive independence samplers, 
discussed in \cite{Giordani10,Holden09,Keith08},
use ideas similar to ISA in the context of MCMC\slash rejection sampling.
The ``cross-entropy method'' also adapts proposal distributions,
but in a more sophisticated way \cite{Rubinstein04}.
ISA also has connections
with sequential Monte Carlo (see section 3.5 of \cite{DelMoral06}).

\section{Implementation and initialization of ISA}
\label{sec:Implementation}
ISA has been discussed in the literature extensively over
the past decades (see above).
In this paper we focus on strategies for dealing with the pragmatic aspects of applying ISA to ``realistic''
problems in science and engineering
and discussing the implementation of the algorithms,
rather than devising fundamentally new algorithmic ideas.
We are targeting posterior distributions
that are difficult to deal with by the MCMC\slash importance sampling technology reviewed above.
More specifically, we wish to create ISA that 
are efficient under the following conditions.
\begin{enumerate}
\item 
The priors are ``broad'', i.e., using the prior distribution
as proposal distribution leads to the collapse of importance sampling.
\item
We can compute first derivatives of the negative log-posterior distribution $F(\theta)$, 
but not second derivatives.
This situation is often encountered in complex computational models,
e.g., if adjoints are available, but second-order adjoints are not.
\item
The negative log-posterior distribution $F(\theta)$
is not convex and may have multiple minima
that are not well separated (see also figure~\ref{fig:IllustrateAlgorithm}).
\item
The effective dimension, $n_\theta$ of the problem is moderate, 
perhaps tens, rather than hundreds or even thousands,
and that the model $\mathcal{M}$ can be evaluated 
with modest computing resources.
\end{enumerate}
As we will demonstrate, ISA can be effective 
under this conditions,
but other, non-iterative, importance sampling methods
or MCMC may be slow to converge.

\subsection{The basic algorithm}
Our basic ISA, 
summarized in pseudo-code in algorithm~\ref{alg:ISA_general},
is as follows.
We first decide on a class of proposal distributions.
In principle, one could choose different classes of proposal distributions at each step in the iteration,
e.g., a log-normal distribution during the first step, 
and a uniform distribution at the next step and so on.
However, here we consider only the case where the class of 
proposal distributions is fixed
to either Gaussians or multivariate $t$-distributions.
The reason is that both distributions are defined by a mean
and covariance which we can easily compute from the samples.
The utility of multivariate $t$-proposal distributions is discussed in 
a more general setting in \cite{Owen}.
\begin{center}
\begin{algorithm}[t]
\caption{Iterative importance sampling algorithm}
\begin{algorithmic}
 \STATE Given: an unweighted ensemble $\left\{\theta^0_j\right\}$, $j=1,\dots,N_e$
 \STATE Construct a proposal distribution $q^0$ from the ensemble.
 \STATE Set $k=0$, and define a tolerance $\text{tol}$
 \WHILE{no convergence reached}
	\STATE Importance sampling with proposal distribution $q^k$ to get $\left\{\theta_j^{k+1},w^{k+1}_j\right\}$, $j=1,\dots,N_e^k$
	\STATE Compute $R^{k+1}$ by equation~(\ref{eq:R}) for this set of weighted samples
	 \IF{$R^{k+1}/R^{k}<\text{tol}$}
		\STATE Convergence is reached
	\ELSE
		\STATE Construct proposal distribution $q^{k+1}$ from $\left\{\theta_j^{k+1},w^{k+1}_j\right\}$, $j=1,\dots,N_e^k$ (see text)
	\ENDIF
	\STATE Set $k\leftarrow k+1$
 \ENDWHILE
\end{algorithmic}
\label{alg:ISA_general}
\end{algorithm}
\end{center}

Given a set of initial samples,
which we obtain using the initialization strategies discussed below,
we compute the sample mean and covariance,
which define the initial proposal distribution $q^0$.
We then generate $N_e^0$ weighted samples $\{\theta_j^1,w^1_j\}$, $j=1,\dots, N_e^0$,
by importance sampling with the proposal distribution $q^0$,
and compute the (weighted) sample mean and covariance 
to define a new proposal distribution $q^1$,
and repeat.
During the iteration, we monitor the quality measure $R$ in equation~($\ref{eq:R}$),
and stop the iteration when~$\vert R_{k+1}-R_{k}\vert$ drops below a specified tolerance. 
This stopping criterium is similar to the one discussed in \cite{OhBerger92,OhBerger89}.

Note that we cannot guarantee that $R$ decreases monotonically
during the iteration when $N_e^k$ is finite,
in part because we must approximate $R$ in equation~(\ref{eq:R})
by averaging; thus, we only know $R$ up to sampling error,
which can be substantial if $R$ is large and $N_e^k$ is not.
Indeed, a large $R$ is always uncertain,
since an approximation of $R$ by $N_e$ 
samples is always less than $N_e$,
however, the ``true'' $R$ may be much larger.
Since the iteration is stopped when the change in $R$ from one
iteration to the next is small, one must consider sampling error 
when choosing the tolerances.

\subsection{Sample sizes and convergence of ISA}
\label{sec:ISAConvergence}
One must pick an ensemble size $N_e^k$ at each step of the distribution.
We advocate choosing this number as large as possible for two reasons.
The first reason is computational.
Since the iteration of ISA is serial,
we wish to use ISA with a small number of iterations
and large $N_e$ at each iteration to leverage parallelism. 
However, the assumption that a large number of samples can be generated 
``easily'' at each iteration restricts the classes of problems we can solve with ISA.
For example, if a super-computer is needed to evaluate the model $\mathcal{M}$ once,
then ISA as described here is not feasible.
Perhaps it is fair to state that ISA may enable parameter estimation
using massively parallel computers, if the forward problem,
i.e., evaluation of the model $\mathcal{M}$, can be done
using a small fraction (order of 0.1\%) of the available resources.

Our second reason for choosing large $N_e^k$
is motivated by the convergence of ISA.
Suppose the support of the proposal distribution at iteration $k$,
$q^k$, includes the support of the posterior distribution 
(which is a routine assumption).
Then weighted samples from the proposal distribution
$q^k$ generate an ``empirical estimate'' of the proposal distribution
that converges in the sense that
weighted averages converge to expected values with respect to
the posterior distribution as $N_e^k\to\infty$.
This follows from the usual importance sampling theory, see, e.g., \cite{CH13},
even if the proposal distribution depends on past samples (see also \cite{CappeEtAl04}).
The next iteration is guaranteed to converge in the same sense,
and under the same routine assumptions.
Thus, ISA converges within one step,
provided one generates sufficiently many samples,
i.e., if $N_e^k$ is large enough.

The importance of a large $N_e^k$ at each step can
be illustrated further under simplifying assumptions.
Suppose that the posterior distribution is Gaussian with variance $\Sigma$,
and that the proposal distribution at iteration $k+1$ is also Gaussian.
We neglect error in the mean, 
i.e., we assume that the mean of the proposal distribution is exactly equal 
to the mean of the posterior distribution,
which can be justified because small errors in the mean
are less severe than small errors in the covariances \cite{Owen}.
We further assume that the covariance matrix,
computed from the weighted samples from the $k$th iteration,
is $\hat{\Sigma}=(1+\varepsilon)\Sigma$, where $\varepsilon$ is a small number.
Under these assumptions, the weights at iteration $k+1$ are
\begin{equation*}
	w^{k+1}\propto \exp\left(-\varepsilon\,\frac{1}{2}x^t\Sigma^{-1}x\right),\quad x=\sim\mathcal{N}(0,\hat{\Sigma}),
\end{equation*}
which leads to 
\begin{equation*}
	R^{k+1}=\left(\frac{1+\epsilon}{\sqrt{1+2\varepsilon}}\right)^{n_\theta} \approx 1+\varepsilon^2\cdot \left(\frac{n_\theta}{2}\right),
\end{equation*}
where $n_\theta$ is the number of parameters we wish to estimate, i.e., the dimension of $\theta$.
We can interpret $\varepsilon$ as ``sampling'' error, and assume $\varepsilon$ is inversely proportional
to the square root of the effective sample size $N_\text{eff}$ in equation~(\ref{eq:R}).
Thus,
\begin{equation*}
	R^{k+1}\approx 1+R^{k}\;\frac{C\cdot n_\theta}{2N_e^k},
\end{equation*}
where $C$ is a constant.
The above equation indicates that the efficiency of the iteration is directly
controlled by how many samples $N_e^k$ one can generate at each step.
However, the expansions above are only valid if $\varepsilon$ is ``small'', 
which means that sampling error is small and the algorithm 
is already performing well with small $R^k$.
In fact most theoretical considerations rely on large ensemble sizes $N_e^k$,
but in practice the behavior of the algorithm
when the ensemble size is \emph{small}
is more interesting.
We study these practical issues in the context of examples
(see sections~4 and~5 below).

\subsection{Initialization}
Initialization of ISA  has a significant impact on the success of the algorithm.
Indeed, the iteration cannot proceed if the initial proposal distribution $q^0$ yields only
one effective sample, so that the variance that defines $q^1$ is zero,
i.e., the algorithm has ``collapsed'' during the first step. 
We suggest two strategies for initializing ISA.

The first initialization strategy is to run MCMC.
However,
if the MCMC has been run long enough to generates sufficiently many samples to characterize the posterior, then
the problem has already been solved
and refinement of a proposal distribution by iteration is unnecessary.
We thus consider ``short'' MCMC runs, and in 
the examples and applications below this yielded good results.
The initialization strategy combines the robustness of MCMC
with the computational advantages of importance sampling.
In this context, one can think of ISA as a means to ``fill out''
an under-resolved MCMC computation, 
while making effective use of massively parallel machines.
This can yield dramatic reductions in run times,
from weeks or months of serial MCMC on a workstation,
to minutes on an HPC system.

The second initialization strategy we propose is to construct a Gaussian mixture model (GMM), 
for which each mixture component is centered at a mode of the posterior distribution.
The various modes can be found in parallel,
by initializing a suitable optimization code
at various locations in parameter space obtained by, e.g.,  
sampling the prior distribution.
The covariances of the mixture are inverses of the approximate Hessians
of the negative log posterior distribution $F(\theta)$ at the modes.
The weights of the mixture components are
\begin{equation}
\label{eq:GMM_weights}
 \psi_j = \frac{\exp(-\phi_j)}{\sum_{i=1}^n\exp(-\phi_i)},
\end{equation}
where $n$ is the number of minima we found and $\phi_i$, $i=1,\dots,n$
are the local minima of $F$.
Alternatively, one could also use weights based on
normalizing constants of the Gaussians, 
but we do not pursue this idea further.

One may find that the optimizer,
initialized at two different points in parameter space finds the same minimum.
We distinguish minima by considering the distance
\begin{equation*}
	d_{i,j} =\left(\mu_i-\mu_j\right)^T H^{-1}_i \left(\mu_i-\mu_j\right),
\end{equation*}
where $\mu_k$ are the various minima we find,
and where $H_k$ are the corresponding approximate Hessians
with elements $(H_k)_{ij} = \partial^2 F/\partial\theta_i\partial\theta_j\vert_{\mu_k}$,
$i,j=1,\dots,n$.
We say that a minimum of $F$ at $\mu_i$ is different from the one
at $\mu_j$ if $d_{i,j}$ is more than an application-specific threshold.
Other, more sophisticated assessments may be needed in other applications
to account for asymmetries, since
\begin{equation*}
	\left(\mu_i-\mu_j\right)^t H^{-1}_i \left(\mu_i-\mu_j\right) \neq \left(\mu_i-\mu_j\right)^t H^{-1}_j \left(\mu_i-\mu_j\right),
\end{equation*}
However, in the applications we consider below, our simple strategy is sufficient.

\subsection{Illustration by a 2D toy problem}
\label{sec:toy}
\begin{figure}[tb]
\centering
\includegraphics[width=.7\textwidth]{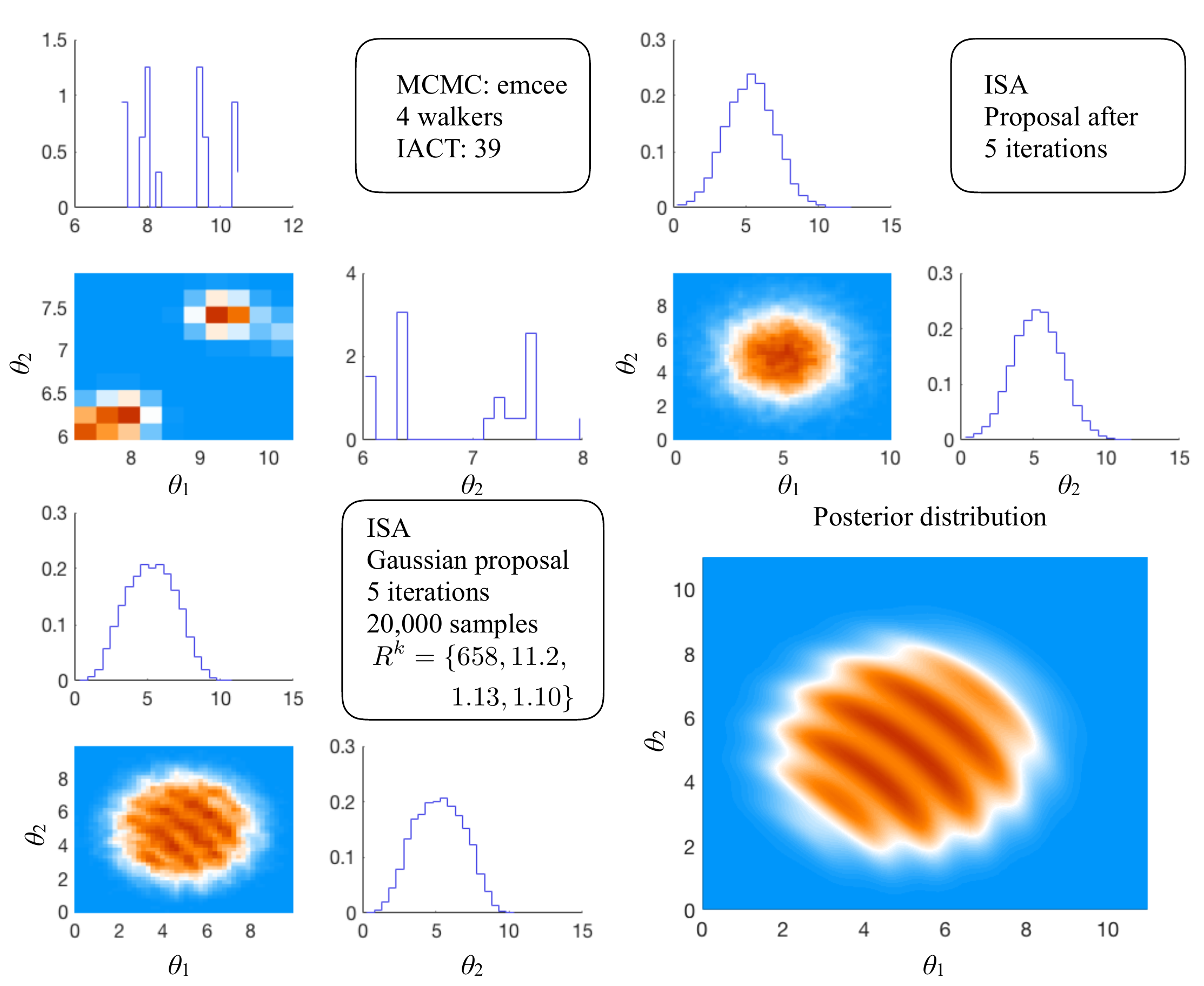}
\caption{
Illustration of ISA initialized by a short MCMC run.
Top row, left: triangle plot of the first 20 MCMC samples.
Top row, right: triangle plot of the Gaussian proposal distribution obtained after 5 iterations.
Bottom row, left: triangle plot of the posterior distribution obtained by ISA.
Bottom row, right: the posterior distribution, bird view.
}
\label{fig:ISA_Illustration_Hammer}
\end{figure}

We illustrate ISA by a simple ``toy'' problem,
where we estimate a two-dimensional parameter vector $\theta$.
The prior is uniform on the cube $\left[0,11\right] \times \left[0,11\right]$,
and the likelihood is such that the posterior is given by 
\begin{equation*}
	p(\theta\vert z) \propto\exp\left(-	F(\theta)\right), \quad \text{where } 
	F(\theta)= 10^{-2} \,\vert\vert \theta-[5; 5] \vert\vert^4+0.2\, \sin(5\vert\vert \theta\vert\vert),
\end{equation*}
and where $\vert\vert x\vert\vert=\sqrt{x^Tx}$ is the 2-norm.
The posterior distribution is illustrated in the bottom right panel of Figure~\ref{fig:ISA_Illustration_Hammer}.
Note that there are several modes that carry significant probability mass,
but the modes are not well separated.

We apply ISA with Gaussian proposal distributions
to obtain samples from this posterior distribution.
We initialize the iteration by samples we obtain
by the emcee implementation \cite{Hammer}
of the affine invariant MCMC ensemble sampler \cite{GoodmanWeare10}.
We use four walkers, initialized by drawing samples from the uniform prior distribution. 
The first 20 samples are used to initialize ISA.
A ``triangle plot'' of these 20 samples, 
which shows histograms of all
one-dimensional marginals and the two-dimensional PDF,
of a given set of samples, is shown in the upper left panel of Figure~\ref{fig:ISA_Illustration_Hammer}.
We perform 5 iterations, each with $N_e=20,000$ samples.
The quality measure $R$ converges quickly to about $R\approx 1.1$,
by the sequence  $\{R^k\}=\{658,11.2,1.13,1.10\}$.
A triangle plot from the Gaussian proposal we obtained after five iterations is shown in
the upper right panel of Figure~\ref{fig:ISA_Illustration_Hammer}.
A triangle plot of 20,000 weighted samples of this Gaussian proposal 
are shown on the lower left panel of the figure,
and one can clearly identify multiple modes in the posterior distribution.
Moreover, the efficiency of importance sampling is almost perfect in this example:
since $R\approx 1$, almost all samples are effective samples.
Emcee hammer on the other hand is characterized by an integrated 
auto correlation time of about 39, which means that only one in 39 samples is effective.
We estimated the integrated auto correlation time from an emcee run
with four walkers, each taking 100,000 steps.

\begin{figure}[tb]
\centering
\includegraphics[width=.45\textwidth]{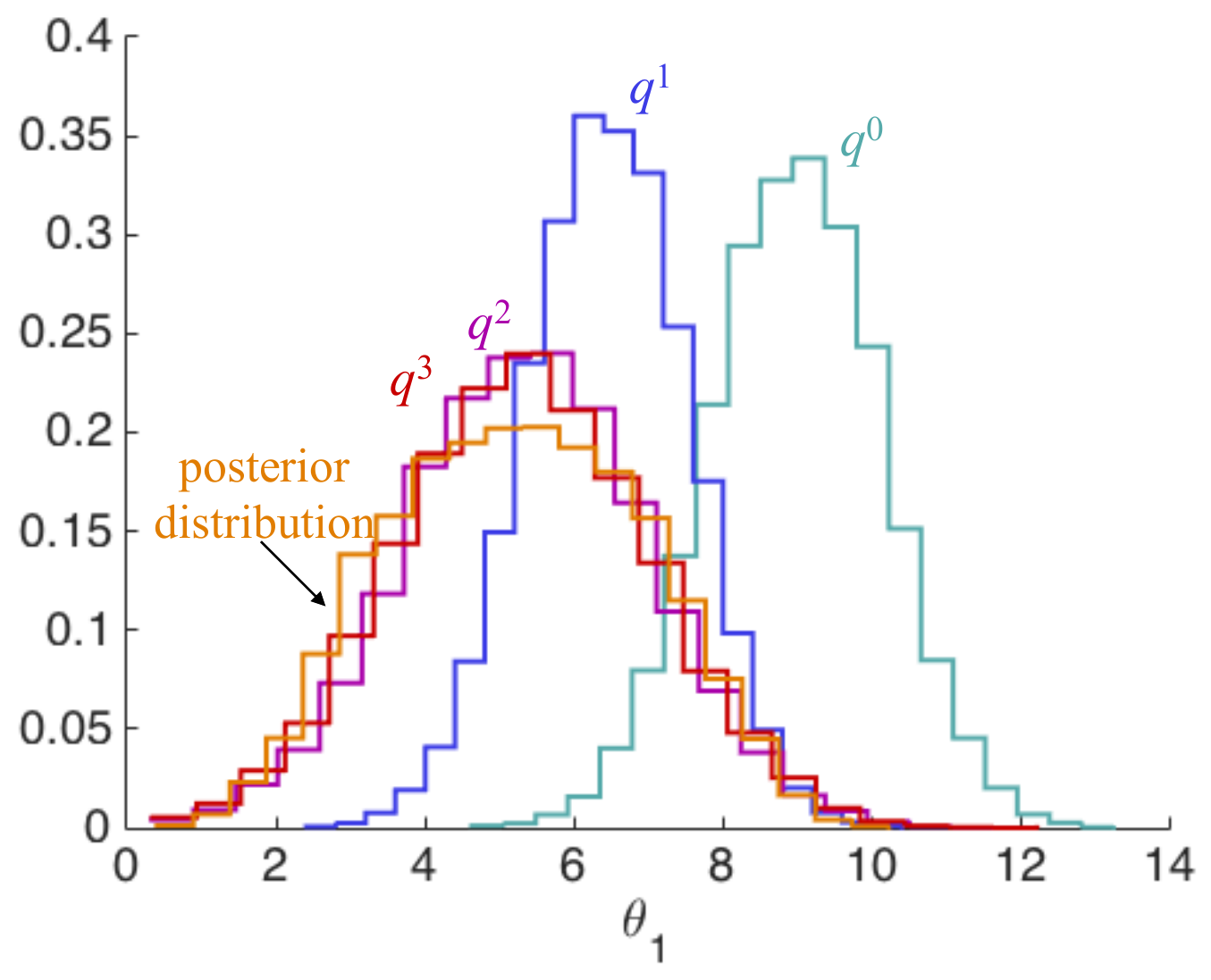}
\caption{Convergence of the proposal distributions of ISA.}
\label{fig:ISA_convergence}
\end{figure}

The convergence of ISA is further illustrated in Figure~\ref{fig:ISA_convergence},
where we show histograms of one-dimensional marginals of the proposal distributions $q^0,q^1,q^2,q^3$,
obtained during the iteration, as well as a histogram of the corresponding marginal of the posterior distribution.
We observe that the high probability region of the  initial proposal distribution is far from
that of the posterior distribution. 
However, ISA quickly corrects this during iteration
and moves the Gaussian proposal distributions closer to the posterior distribution.

We also tested ISA initialized by a GMM,
which we obtain by performing 100 minimizations of $F$, 
each initialized by a sample from the uniform prior density.
We then compare the minima by the distance $d_{i,j}$
as explained above, and find five distinct minima of $F$.
At each distinct minimum, we approximate the Hessian by finite differences
(we chose computationally more 
efficient methods to approximate Hessians in the applications below).
The distinct minimizers and corresponding Hessians define the GMM
with five mixture components,
where the weights of each component is given by~(\ref{eq:GMM_weights}).
We use this GMM as the proposal distribution and compute 20,000 samples to obtain $R\approx 2.59$.
Using the GMM mixture distribution is thus less effective than using a few iterations of ISA as above.
Nonetheless, we can use a few samples of the GMM to initialize ISA.
Here we use the 50 samples from the GMM to initialize ISA
and then iterate using Gaussian proposal distributions with 20,000 samples at each step.
The Gaussian proposal distribution we obtain after three iterations 
is characterized by a quality measure $R = 1.10$ 
(the sequence is $R=\{3.27, 1.1, 1.1\}$). 
The experiments indicate the iterations of ISA yield 
similar proposal distributions after a few iterations
for both initialization strategies.

\newpage
\section{Application 1: subsurface flow}
\label{sec:Subsurface}
\subsection{Background and problem set-up}
Subsurface flow models are extensively used in hydrological and geosciences applications, such as the prediction of contaminant fate during a chemical spill, and the estimation of oil production from an oil recovery process.
Robust predictions of quantities of interest require accurate characterization of the model parameters, such as the material properties and process-related parameters. 
Difficulties in direct measurements of these parameters necessitate the use of inverse modeling tools to infer them based on sparse measurements of observables that can be modeled.  

\begin{figure}[tb]
\centering
\includegraphics[width=.8\textwidth]{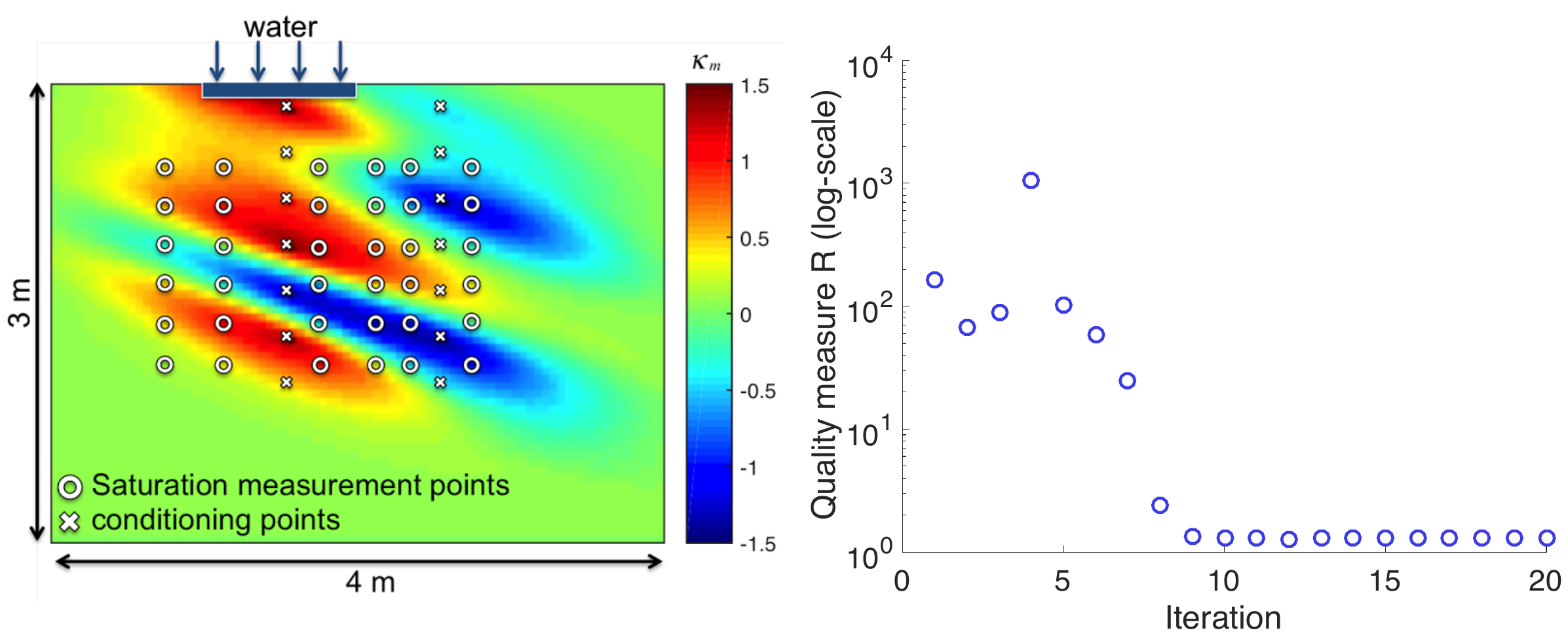}
\caption{Left: Computational domain and locations of the conditioning points for Kriging; 
the color map illustrates the distribution of the permeability modifier 
when the values at the pilot points are given by  
$\theta_\text{ref}$.
Right: quality measure $R$ during ISA iteration.}
\label{fig:subsurface_setup}
\end{figure}

In this section, we use ISA in a synthetic subsurface flow problem
to estimate the permeability distribution, and the associated uncertainty, in the vadose zone from sparse measurements of
saturations and changes in water volume in a domain that covers a two-dimensional 4m$\times$3m (width$\times$depth) spatial domain and is discretized by a mesh with uniform grid cells of size $0.05 \times 0.05$m. Water is released from a pond into a heterogeneous vadose zone that has a water table at depth 3m. 
The saturation of the heterogeneous soil is initially in capillary-gravity equilibrium. 
A one-day ponded infiltration period is followed by a one-day water redistribution period. 
During infiltration, the water level in the pond is maintained at 2cm,
and saturations at 36 locations, shown in the left panel of Figure~\ref{fig:subsurface_setup},
are measured every two hours (25 snapshots) during the two days of infiltration and redistribution. 
The total volume of water flowing out of the pond is also measured at the same frequency, yielding a total of 925 observations that we combine into a vector $z$ as in equation~(\ref{eq:ModelToData}). 

Our goal is to determine the permeability distribution of the vadose zone based on the above measurements.  
We describe the permeability 
by the pilot-point method~\cite{ramarao95,gomez97} with fourteen pilot points,
located in two columns near the pond,
as shown in the left panel of Figure~\ref{fig:subsurface_setup}.
The pilot points serve as adjustable conditioning points for determining the permeability modifier ($\kappa_m$) in each grid cell through a Kriging process that utilizes a spherical semivariogram model with fixed parameters~\cite{deutsch92}. The permeability in each grid cell is then given by $\kappa_{\rm ref} \times 10^{\kappa_m}$,
where $\kappa_{\rm ref} = 0.316 \times 10^{-11}$ m$^2$ is a reference permeability value.
The inverse problem we consider amounts to
estimating the parameter vector $\theta =\{\kappa_{m,i},1\le i\le 14\}$, where $\kappa_{m,i}$ are the $\kappa_m$ values at the pilot point locations, based on two days worth of saturation and flow data.

We assume a broad Gaussian prior with zero mean and 
a diagonal covariance matrix $C_{pp}$ where all 14 diagonal elements are equal to nine.
The likelihood is defined by equation~(\ref{eq:ModelToData}),
where the model $\mathcal{M}$ is based on Richards' equation~\cite{Richards31}.
Numerically, we solve this equation using the finite-volume simulator TOUGH2~\cite{Pruess99} 
with relative permeability and capillary pressure functions as in~\cite{vGenuchten80}.
The errors $v$ in equation~(\ref{eq:ModelToData}) are
Gaussian with mean zero and a diagonal covariance $C_{zz}$.
The standard deviations that define the covariance are assumed to be
0.01 for the saturation measurements,
and 0.005 m${^3}$ for the water volume measurements.
We perform experiments with synthetic data,
i.e., we assume a reference permeability field defined by
\begin{equation*}
	\theta_\text{ref} = 
	\{1.5, .5, 1., 1.5, -1.5, 1.5, 0.5, -0.5, -0.5, -1.5, -0.001, 0.5, -1.5, -0.5\},
\end{equation*}
and then perform a simulation of the above infiltration. The reference set of parameters defines the permeability modifier field shown in 
the left panel of Figure~\ref{fig:subsurface_setup}. The synthetic observation vector $z$ is then given by $\mathcal{M}(\theta_{\rm ref}) + \epsilon$ where $\epsilon$ is an instance of $v$ in equation~(\ref{eq:ModelToData}).

\subsection{Results}
We initialize ISA by a single Gaussian.
We define this Gaussian by solving the nonlinear least squares problem
$\min_\theta F(\theta)$, 
where $F$ is the negative logarithm of the posterior (up to an additive constant)
\begin{equation}
F(\theta)=\frac{1}{2}(z-\mathcal{M}(\theta ))^TC_{zz}^{-1}(z-\mathcal{M}(\theta ))+\frac{1}{2}\theta^TC_{pp}^{-1}\theta. \label{eqn:F_subsurface}
\end{equation}
Here $\mathcal{M}(\theta)$ is the output of the TOUGH2 simulator. 
The minimizer of $F$ defines the mean,
and the approximate Hessian, evaluated at the minimum,
defines the inverse of the covariance matrix.
We approximate the Hessian by neglecting contributions from second derivatives,
i.e., we set the covariance equal to the inverse of
 $2\,J^TJ$, where $J$ is the Jacobian matrix of the nonlinear least squares problem.
Numerically, we rely on Levenberg-Marquardt \cite{levenberg44,marquardt63}  (lsqnonlin function in MATLAB)
to solve the optimization problem and use finite differences for computing the first derivatives. 
For this particular problem, 20 iterations are needed to find a minimum. 
We note that the evaluations of 
the negative log-posterior~(\ref{eqn:F_subsurface}) and 
of the Jacobian $J$ in each iteration can be easily parallelized. 

We start ISA using the Gaussian proposal distribution and generate 5,120 samples at each iteration.
The likelihoods and weights of the samples are evaluated in parallel
using 160 nodes, each with 32 cores, on Cori, a supercomputer at NERSC. 
With this set-up, the quality measure $R$ converges from 163 to 1.3 in 9 iterations,
as illustrated by the right panel of Figure~\ref{fig:subsurface_setup}. 

We compare ISA to solving the same problem by the affine invariant ensemble sampler
\cite{GoodmanWeare10}, as implemented in the emcee code~\cite{Hammer}.
We use a parallel implementation of the algorithm with 126 walkers,
which allows us to use 2 nodes on Cori efficiently.
We discard the first 2,500 steps as ``burn-in''
and take another 2,500 steps per walker. 
The computational cost thus amounts to a total of  630,000 likelihood evaluations, 
each requiring an evaluation of the TOUGH2 model that takes approximately 3 minutes.  
The integrated auto-correlation times for $\theta$ vary from 84 to 532, 
indicating that we have 3,750 to 592 ``effective" samples.  
This should be compared to the quality measure of $R\approx 1.3$ after nine iterations of ISA.
This means, in particular, that after nine iterations, 
almost every Gaussian sample contributes effectively
to the numerical approximation of the posterior density.
To get there, about 46,080 samples are used,
and these can be evaluated in parallel in nine batches of 5,120 samples. 
Our results demonstrate that ISA can solve this parameter estimation problem effectively, 
while achieving greater level of parallelism than MCMC,
which reduces wall-clock time. 
Including the cost of the initial minimization of (\ref{eqn:F_subsurface}), 
solving the problem by ISA takes about 90 minutes, 
while the emcee code requires about 20 days.

\begin{figure}[t]
\centering
\includegraphics[width=.55\textwidth]{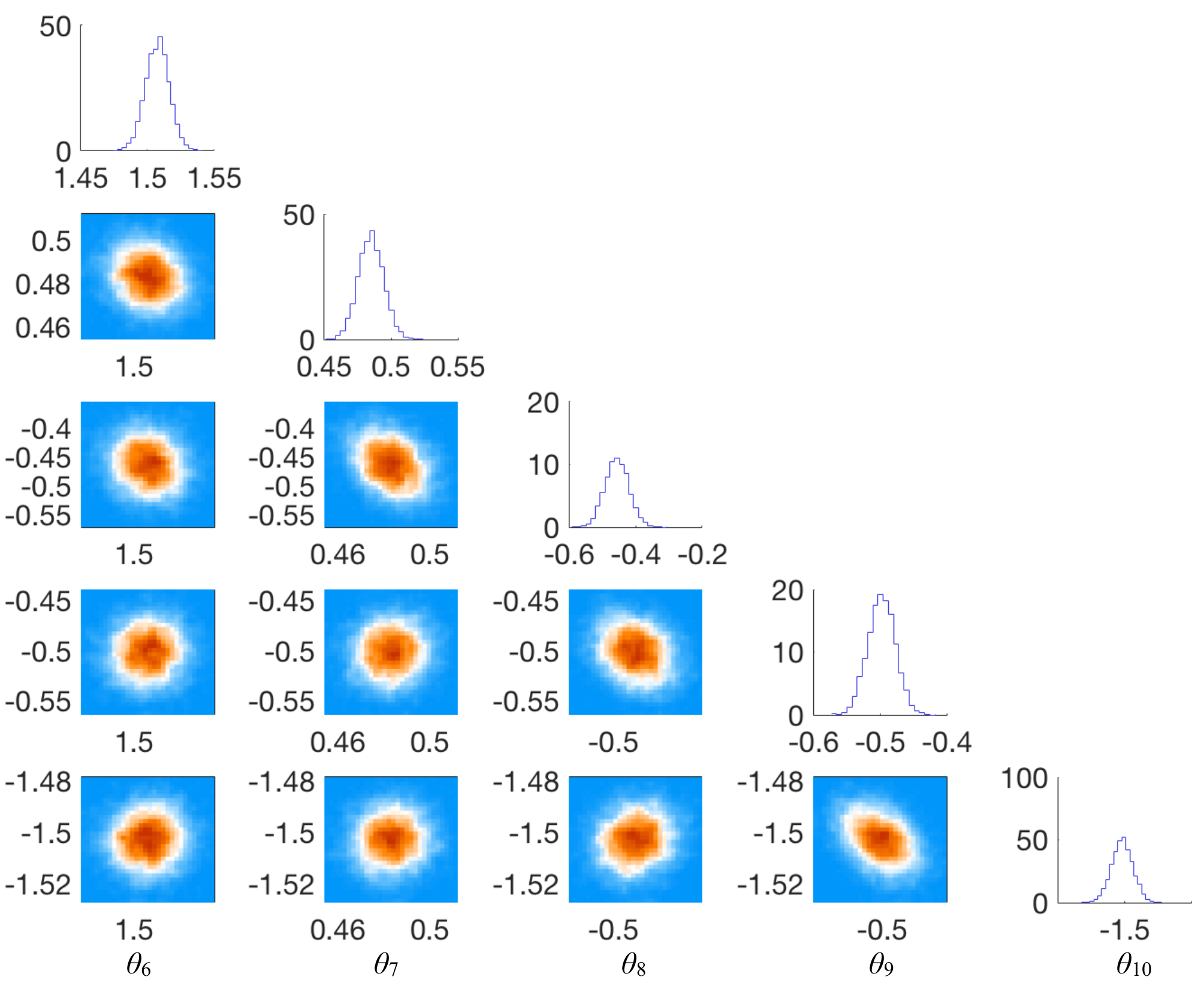}
\caption{Triangle plot of posterior distribution obtained by ISA after nine iterations, each with 5,120 samples. 
The plot shows histograms of the parameters 
$\{\theta_6,\theta_7,\theta_8,\theta_9,\theta_{10}\}$ (diagonal),
as well as two-dimensional histograms of all pairs of 
these five parameters (lower triangle).}
\label{fig:subsurface_corner}
\end{figure}

Note that in the context of ISA,
a quality measure $R$ close to 1 can
indicate that the posterior is nearly Gaussian.
The reason is that we use a Gaussian proposal distribution,
and $R\approx 1$ implies that the posterior distribution
is also nearly Gaussian. 
This is illustrated by the triangle plot in Figure~\ref{fig:subsurface_corner}.
The plot shows histograms of the five parameters $\theta_6$--$\theta_{10}$,
as well as two-dimensional histograms of pairwise 
combinations of these five parameters. 
Histograms of those parameters that are not shown in the figure are qualitatively similar.
The figure suggests that the variables are indeed nearly Gaussian. 
While we can conclude, after the fact, that
other methods specifically designed for nearly Gaussian problems
(see literature review) may perform better on the present problem, 
the distributions of the parameters are not typically known~{\em a~priori}. 
On the other hand, 
the Gaussian proposal distribution whose mean is the posterior mode
and whose covariance matrix is derived from the Gauss-Newton solution
if the optimization problem, exhibits an $R$-value much larger
than what we find after a few iterations of ISA.

\section{Application 2: combustion modeling}
\label{sec:Combustion}
\subsection{Background}
Predictive combustion modeling entails solving a set of fluid conservation equations that
are supplemented with models for thermodynamics, transport and chemical kinetics.  
For this study we assume an ideal gas mixture, 
gradient-driven transport,
and a set of reversible Arrhenius rate expressions
for chemical reactions.
Here, we focus on chemical combustion on hydrogen combustion.
Compared to fuels involving complex
hydrocarbons or biofuels, the hydrogen-oxygen case is composed of
a relatively small number of chemical species and reactions.
For that reason, and perhaps since hydrogen-oxygen
represents an important building block for more complex fuels,
the basic form of the numerical model for hydrogen-oxygen systems, 
including the set of specific fundamental chemical reactions, 
is well-established.
The parameters that define the model, however, remain open to considerable debate,
as indicated by the number of models under development
(e.g., \cite{DavisEtal2005,GlarborgAlzueta1998,SmithGolden,LiDryer2004,Miller,OConnaire2004}).
Note that the various parameter values reported over the years
for each line of model development result from a long history of incremental optimization
efforts that incorporate emerging experimental data sets.  
As a result, while different parameter sets may lead to similar predictive capabilities
within a common range of their respective validated conditions, the models they define
can generate dramatically different predictions outside that range (see the discussion,
e.g., in \cite{Grcar2008}, particularly related to Fig.~1, and in \cite{DongEtal2005},
related to Fig.~4).

The numerical model used here consists of a systems of differential equations
whose solution requires a database of parameters that
characterize the relevant thermodynamic relationships, 
transport coefficients, and kinetic rate coefficients for the component gases.
Our goal is to use ISA to incorporate measured experimental data into 
characterization of the uncertainty of a selected set of nine parameters that
describe important behavior of this particular numerical model.

To a large extent, our study follows a recent effort to optimize 28 parameters of
a syngas combustion model \cite{DavisEtal2005}.  In that work, the
model was tuned specifically to predict measured values of 36 experiments,
including constant-pressure flow reactors, constant-volume shock-tubes
and steady, unstrained premixed flames.
In the reference, a solution-mapping technique was employed to approximate
the response of the model using a second-order polynomial.  For each
experiment, an {\it {a priori}} sensitivity analysis was used to
reduce the set of active parameters in the response surface model.  An
objective function was constructed and minimized based on the square
of model prediction errors, normalized by the uncertainty of the
respective measurements.  
Prior knowledge of the parameters was used only to enforce bounds on the variation of
each active parameter, and to inform the construction of the underlying polynomial response surface model.
Implicitly, this work also incorporated prior knowledge in the
choice of what reactions to include in the kinetic model and the choice of which model parameters to hold
fixed during the optimization.
In a Bayesian framework, the referenced optimization approach can be interpreted
loosely as one based on a uniform prior
over a cube whose sides are defined by prior-informed parameter bounds.
In the present work, discussed below, we estimate a subset of the parameters considered in \cite{DavisEtal2005}
and quantify the respective uncertainties using ISA, based on the complete numerical combustion model,
rather than solution-mapping.  However, we make use of many of the same target experiments in
order to quantify our likelihood estimate.

\subsection{Prior and forward model, likelihood and data}
We use the optimized mechanism of \cite{DavisEtal2005},
shown in detail in Table~\ref{tab:srcmech},
to define the mean of a truncated Gaussian prior distribution for the
parameters shown in red.
Since we wish to quantify the extent to which experimental data (see below)
constrains the active parameters, we take the standard deviations of each
parameter to be of the order of its mean value (as shown in the table).
\begin{table}[tb]
\scriptsize
  \centering
  \begin{tabular}{r @{\hspace{1em}} r c l @{\hspace{1em}} r @{\hspace{1em}} r @{\hspace{1em}} r @{\hspace{2em}} r @{\hspace{1em}} r}
\multicolumn{1}{c}{ID}  &
\multicolumn{3}{c}{Reaction} &
\multicolumn{1}{c}{A$^{*}$} &
\multicolumn{1}{c}{$\beta$} &
\multicolumn{1}{c}{$E_a$} &
$\sigma$&
Bounds (L:U)\\ \hline
1   & H + O$_2$       &=& O + OH             & \Red{2.644(16)} & -0.6707 &17041 & 2(16) & 1(16):5(16) \\
2   & O + H$_2$       &=& H + OH             & \Red{4.589(4)}  &  2.7    & 6260 & 1(4)  & 1(4):1(5) \\
3   & OH + H$_2$      &=& H + H$_2$O         & \Red{1.734(8)}  &  1.51   & 3430 & 2(8)  & 4(7):4(8) \\
4   & 2 OH            &=& O + H$_2$O         & 3.973(4)        &  2.4    &-2110 & & \\
5   & 2 H + M         &=& H$_2$ + M          & 1.780(18)       & -1      &    0 & & \\
6   & H + OH + M      &=& H$_2$O + M         & 4.400(22)       & -2      &    0 & & \\
7   & O + H + M       &=& OH + M             & 9.428(18)       & -1      &    0 & & \\
8   & 2 O + M         &=& O$_2$ + M          & 1.200(17)       & -1      &    0 & & \\
9   & H + O$_2$ (+M)  &=& HO$_2$ (+M)        &  &  & & & \\
    &  \multicolumn{3}{@{\hspace{1.2em}} l}{{\em high pressure, K$_{f\infty}$}} & \Red{5.116(12)} & 0.44 & 0  & 5.116(12) & 0:1.6(13)\\
    &  \multicolumn{3}{@{\hspace{1.2em}} l}{{\em low pressure, K$_{f0}$}}      & \Red{6.328(19)} & -1.4 & 0  & 6.328(19) & 0:2(20)\\
    &  \multicolumn{6}{@{\hspace{1.2em}} l}{{\em TROE:} $F_c = 0.5$} & & \\
    &  \multicolumn{6}{@{\hspace{1.2em}} l}{{\em Third-body:}  O$_2$(0.85), H$_2$O(11.89), Ar(0.40)} & & \\
    &  \multicolumn{6}{@{\hspace{1.2em}} l}{\hspace{5em}  He(0.46), H$_2$(0.75)} & & \\
10  & H$_2$ + O$_2$   &=& HO$_2$ + H         & \Red{5.916(5)} &  2.433  &53502 & 5.916(5)& 0:2(6)\\
11  & 2 OH (+M)       &=& H$_2$O$_2$ (+M)    &  & & & & \\
    &  \multicolumn{3}{@{\hspace{1.2em}} l}{{\em high pressure, K$_{f\infty}$}} &  1.110(14) & -.37   & 0 & & \\
    &  \multicolumn{3}{@{\hspace{1.2em}} l}{{\em low pressure, K$_{f0}$}} &  2.010(17) & -.584   & -2293 & & \\
    &  \multicolumn{6}{@{\hspace{1.2em}} l}{{\em TROE:}
              $F_c=0.2654\exp{(-T/94)}$}& & \\
    &  \multicolumn{6}{@{\hspace{1.2em}} l}{\hspace{5em}
              $+0.7346\exp{(-T/1756)}+\exp{(-5182/T)}$}& & \\
    &  \multicolumn{6}{@{\hspace{1.2em}} l}{{\em Third-body:}  H$_2$(6), H$_2$O(6), Ar(0.7), He(0.7)}& & \\
12  & HO$_2$ + H      &=& O + H$_2$O         & 3.970(12)       &  0      &  671 & & \\
13  & HO$_2$ + H      &=& 2 OH               & \Red{7.485(13)} &  0      &  295 & 7.485(13) & 0:2(14) \\
14  & HO$_2$ + O      &=& OH + O$_2$         & \Red{4(13)}     &  0      &    0 & 4(13)     & 0:2(14) \\
15a & HO$_2$ + OH     &=& O$_2$ + H$_2$O     & \Red{2.375(13)} &  0      & -500 & 2.375(13) & 0:1(14) \\
15b & HO$_2$ + OH     &=& O$_2$ + H$_2$O     & 1.000(16)       &  0      &17330 & & \\
16a & 2 HO$_2$        &=& O$_2$ + H$_2$O$_2$ & 1.300(11)       &  0      &-1630 & & \\
16b & 2 HO$_2$        &=& O$_2$ + H$_2$O$_2$ & 3.658(14)       &  0      &12000 & & \\
17  & H$_2$O$_2$ + H  &=& HO$_2$ + H$_2$     & 6.050(6)        &  20     & 5200 & & \\
18  & H$_2$O$_2$ + H  &=& OH + H$_2$O        & 2.410(13)       &  0      & 3970 & & \\
19  & H$_2$O$_2$ + O  &=& OH + HO$_2$        & 9.630(6)        &  2      & 3970 & & \\
20a & H$_2$O$_2$ + OH &=& HO$_2$ + H$_2$O    & 2.000(12)       &  0      &  427 & & \\
20b & H$_2$O$_2$ + OH &=& HO$_2$ + H$_2$O    & 2.670(41)       & -7      &37600 & & \\
  \end{tabular}
  \caption{Arrhenius rate parameters for H$_2$-O$_2$ combustion model (kinetics, and accompanying 
           thermodynamics and transport parameter database taken from \cite{DavisEtal2005}).
           Parameters in red are active for the present study.
           The forward rate constant, $K_f = AT^{\beta}\exp{(-E_a/RT)}$. $^*$The number
           in parentheses is the exponent of 10, i.e., 2.65(16) = 2.65 $\times$ 10$^{16}$.
           For the active parameters, standard deviation of the prior, $\sigma$ and lower
           and upper bounds of the truncated Gaussian are provided in the last two columns.}
  \label{tab:srcmech}
\end{table}
In this way, we generate a broad Gaussian distribution that extends well beyond the upper and lower bounds specified for
each parameter (shown in the rightmost column of the table).
We truncate the Gaussians at the bounds to obtain our truncated Gaussian prior distribution. We
note that these bounds represent are much broader that those considered in \cite{DavisEtal2005}.

Our goal is to sharpen the prior knowledge by incorporating a limited set of experimental data.
We consider the six flow reactors and eight ignition delay experiments labeled
({\em ign1a, ign1b, ign2a, ign2b, ign3a, ign4a, ign4b, ign5a, flw1a, flw2a, flw3a, flw4a, flw8a, flw8b})
in Table~1 of \cite{DavisEtal2005}. The data include Gaussian error estimates,
taken as the standard deviations,
$\sigma$, reported in Table~1 of \cite{DavisEtal2005}.
We assume that the errors are uncorrelated,
i.e., the covariance matrix that defines the errors $v$ in equation~(\ref{eq:ModelToData}) is diagonal.
The data, assembled into a $14\times 1$ vector $z$, and the Gaussian error model are used along with our numerical model 
to define a likelihood by equation~(\ref{eq:ModelToData}).

In equation~(\ref{eq:ModelToData}),
each element of the vector $z$ corresponds to one of the 14 experiments. 
We compute $z$ using point-reactors
subject to constant pressure (flow-reactor) or constant volume (shock-tube) constraints. 
VODE \cite{vode}, is used to integrate either set of equations over time, and the 
solution is monitored until the desired diagnostic is obtained.  
Initial data for the simulations is provided in the description of the experiment in the reference,
as is the event defining the measured quantity. 
For example, experiment {\em flw1a}
reports the mean rate of change of H$_2$ mole fraction during the time interval when it is between
40 and 60 percent of its initial value. 
The computed and measured data, along with the sampled
and prior mean values of the parameters, and the corresponding prior standard deviations and 
experimental measurement errors are combined for each sample vector to form the resulting likelihood.

Our combustion model is accessorized with several options for efficiency and robustness so 
that it can be driven by a number of sampling strategies.  
In particular, the model is allowed to fail gracefully if any component of a sample 
(a vector of nine parameters) is generated outside the bounds specified in Table~\ref{tab:srcmech}.
In this case, the likelihood is treated formally as zero.  Similarly, some combinations of 
parameter values lead to predicted values that are far from the corresponding measured 
data.  Ignition delays, for example, may be prohibitively long to compute efficiently, or too 
short to accurately capture when scanning the solution at a ``reasonable'' sampling frequency.
In either case, the diagnostic analysis routines can identify the failure mode and gracefully
exit, signaling the driver code that the corresponding likelihood is to be taken as zero.
Extra controls are provided/required in the specification of each experiment to help detect
such scenarios.  These auxiliary parameters are adjusted to have minimal effect on the presented data,
but enable robust performance of the model evaluation over millions of sample sets.
For each sample, the list of simulated experiments are evaluated in parallel over compute threads, 
using an OpenMP task manager.  
When a large set of independent sample vectors can be generated at once, 
the likelihood function itself can be evaluated in parallel over samples using MPI, 
where each MPI rank maintains its own version of the model parameter set and evaluates 
its likelihood using a local subset of the available cores.

\subsection{Results}
\begin{figure}[tb]
\centering
\includegraphics[width=1\textwidth]{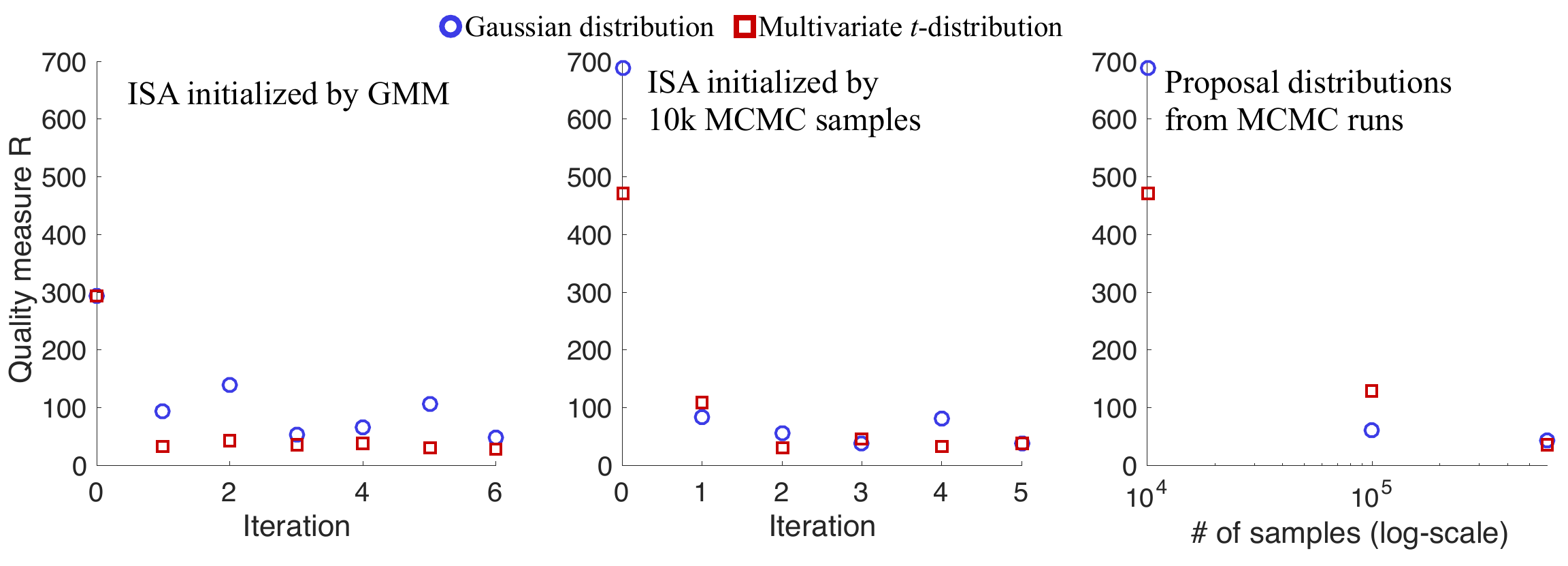}
\caption{
Iterations of ISA.
Left: Quality measure $R$ as a function of iteration number
for ISA initialized by a GMM, using Gaussian or multivariate $t$-distributions 
and $10^5$ samples at each iteration.
Center: Quality measure $R$ as a function of iteration number
for ISA initialized by $10^4$ samples of an MCMC run,
using Gaussian or multivariate $t$-distributions 
and $10^5$ samples at each iteration
Right: Generating Gaussian\slash multivariate $t$-proposal distributions from 
$10^4$, $10^5$, and $6\cdot 10^5$ samples.
}
\label{fig:CombustionIterations}
\end{figure}

As explained in section~\ref{sec:ParameterEstimationProblem},
the prior distribution and a likelihood jointly define a posterior distribution,
and we apply several variants of ISA algorithms to draw samples from this posterior distribution.
We first test initialization by a GMM. 
To construct the GMM we run $10^4$
minimizations of the negative logarithm of the posterior.
The starting point for the optimizations are drawn from the broad prior.
We use the nonlinear least squares code
in the C implementation of ``minpack''  to carry out the required optimizations.
One could also use derivative free optimization,
but we decided to use minpack and compute
derivatives by finite differences.
We use $2,400$ cores on NSERC's supercomputer ``Edison''
to perform the optimization and of our $10^4$ optimization attempts, 
$2,482$ were ``successful''.
The large number of failures is caused our broad prior: 
many of the starting parameter sets we tried were outside of the range
the numerical model can simulate,
which made the optimizer return a failure. 
We identified $12$ minima to be distinct 
by using the criterion outlined in section~\ref{sec:Implementation}.
Here we use a threshold based on a $\chi^2$ distribution
so that the probability of the minima being within the threshold exceeds 95\%.
The resulting Gaussian mixture consists of 12 components,
centered at the minima.
The covariances are computed from approximate Hessians at these minima,
using the same Gauss-Newton approximation as in section~\ref{sec:Subsurface}.
We generate $10^5$ samples from this GMM and evaluate
the corresponding likelihoods using 3,072 cores on Edison.
Evaluating and weighting $10^5$ samples  takes less than 15 minutes of wall-clock time,
and results in the quality measure, $R=293$.
From the weighted samples we generate Gaussian or multivariate $t$-distributions
with shape parameter $\nu=3$ and iterate as described above. 
At each step of the iteration we ``inflate'' the sample covariance by a factor of two
to mitigate sampling error and the underestimation of covariances,
as is common in numerical weather prediction \cite{Anderson2009}.
The inflation parameter could be tuned further,
however we did not pursue this because we obtained good results.
The left panel of Figure~\ref{fig:CombustionIterations}
shows the quality measure $R$ as a function of the iteration number when using a multivariate t-distribution.
The figure suggests that 
ISA produces an effective proposal distribution after the first iteration 
with an associated $R$-value of less than 42.
The precise sequence during iteration is $R=\{ 293, 32, 42, 35, 38, 29, 27\}$. 

Also shown in the left panel of Figure~\ref{fig:CombustionIterations} are results
we obtain by Gaussian proposal distributions at each step of ISA.
We note that the values of $R$ fluctuate more than when using multivariate $t$-distributions,
and the lowest value of $R$ we obtain is also larger than the one we find by using 
multivariate $t$-proposal distributions.
The precise sequence is $R=\{ 293, 94, 140, 54, 66, 105, 47\}$.
Nonetheless, the values of $R$ we obtain,
in combination with the massive parallelism we can leverage at each step of the iteration,
make ISA with Gaussian proposal distributions effective for this problem.

In the center panel of  Figure~\ref{fig:CombustionIterations}, 
we show results we obtain by ISA initialized by 
multivariate $t$ or Gaussian distributions computed from $10^4$ samples of an MCMC run.
As before we use the emcee implementation of the affine invariant ensemble sampler 
\cite{Hammer,GoodmanWeare10}.
We describe the details of our MCMC run below,
as it also serves as a reference solution.
In the context of ISA, 
we note that the initial Gaussian\slash multivariate $t$-distribution
yields a significantly larger $R$ than the GMM we constructed above.
The iteration however can quickly reduce $R$ to values of around
$R=42$ (Gaussian proposal), and $R=35$ (multivariate $t$-proposal).
As above, we use $10^5$ samples for each iteration.
Initializing ISA with only a ``few'' MCMC samples
requires at least two iterations to obtain $R$-values below 50,
whereas we could reach low values of $R$ after only one iteration 
when we initialize ISA by a GMM.
Constructing the GMM is computationally more expensive
due to the many optimizations we need to run,
however constructing the GMM can leverage parallelism effectively.
We conclude that, for this problem and computational set-up,
the initialization by GMM is the more effective strategy.

\begin{figure}[tb]
\centering
\includegraphics[width=.7\textwidth]{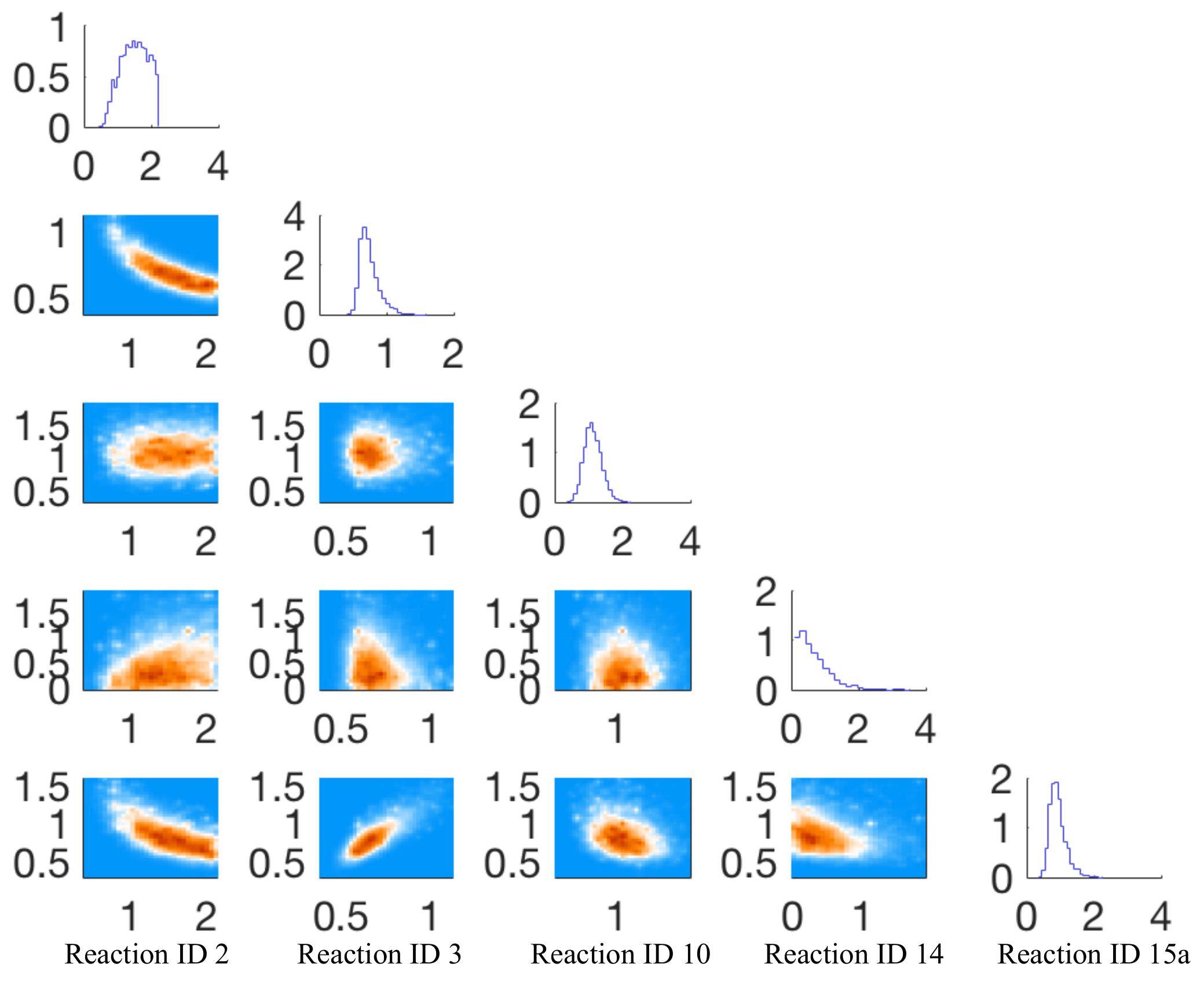}
\caption{
Triangle plot of posterior distribution obtained by ISA after six iterations, each with $10^5$ samples. 
The plot shows histograms of reaction IDs 
$\{2,3,10,14,15a\}$ (diagonal),
as well as two-dimensional marginals of all pairs of 
these five parameters (lower triangle).}
\label{fig:CombustionTriangle}
\end{figure}
We compare the ISA iterations to constructing
Gaussian and multivariate $t$-proposal distributions
from longer MCMC runs.
These tests will
(\emph{i})
illustrate the degree to which ISA improves sampling efficiency 
compared to longer MCMC runs to obtain comparable results;
and (\emph{ii})
demonstrate the convergence of ISA on this problem.
In the right panel of Figure~\ref{fig:CombustionIterations}
we show the values of $R$ we compute for proposal distributions
as a function of the number of MCMC samples we use
to compute the mean and covariances of these distributions.
Specifically, we use $10^4$, $10^5$, and $6\cdot 10^5$
to obtain $R=\{689,60, 42\}$ for Gaussian proposals
and $R=\{471,130, 34\}$ for multivariate $t$-proposals.
We thus obtain similar values of $R$ as when iterating by ISA,
provided that sufficient MCMC samples are used to construct the proposal distribution.
This suggests that ISA can indeed find 
Gaussian or multivariate $t$ approximations of the posterior distribution
that can be used to generate posterior samples with efficiency of $R\approx 40$.

For our MCMC runs, we use a serial implementation of the emcee code,
but parallelize the likelihood evaluations across the various experiments
(flow reactors and ignition delays).
In this way, we could effectively utilize the eight cores on a workstation.
Generating $6\cdot 10^5$ samples then requires about one month.
For comparison, initializing ISA by $10^4$ MCMC samples 
requires about 10 hours of MCMC for generating the initial ensemble,
then 15 mins per iteration using 3,072 cores.
Initializing ISA by GMM further reduces the 
wall-clock time required to generate the initial ensemble.
On the other hand the emcee code also can be parallelized to a limited extent,
which would reduce run-times and improve overall performance.
Nonetheless, it is difficult to achieve the same level of parallelism as ISA
by emcee or other MCMC codes.
Moreover, the ``shape'' of the posterior distribution is difficult
for the emcee code, 
which results in relatively long integrated auto-correlation times of a few hundred samples.
Precise estimates of autocorrelation times are difficult to come by,
however here we are satisfied with the estimates provided by the emcee code
(which differ significantly across the nine variables).
An autocorrelation of a few hundred means
that the ``efficiency'' of the emcee code is about one effective sample
for every few hundred we try.
This should be compared to the quality measure $R\approx 42$ of ISA (after a few iterations),
which suggests that one in about 40 samples is effective.

Finally we show a triangle-plot of the posterior distribution we obtain
by ISA after 6 iterations initialized by a GMM.
The plot shows histograms of five parameters
(Reactions IDs $\{2,3,10,14,15a\}$),
as well as two-dimensional marginals of pairwise 
combinations of these five variables. 
In the plot, the samples are scaled by the prior mean
as shown in Table~\ref{tab:srcmech}.
We note that there is significant probability mass away from the prior mean,
i.e., away from the value one.
Moreover, the posterior mass is concentrated
on narrow ridges in parameter space,
which suggests that posterior variances are reduced significantly compared to the broad prior.
We also note that there are rather strong correlations among some of the variables,
e.g., reaction IDs two and three,
and such correlations are absent in the prior.
Thus, the experimental data we use here do indeed constrain the parameters we wish to estimate,
and ISA can identify the parameter regimes that are consistent with the data we consider.

\begin{figure}[tb]
\centering
\includegraphics[width=.4\textwidth]{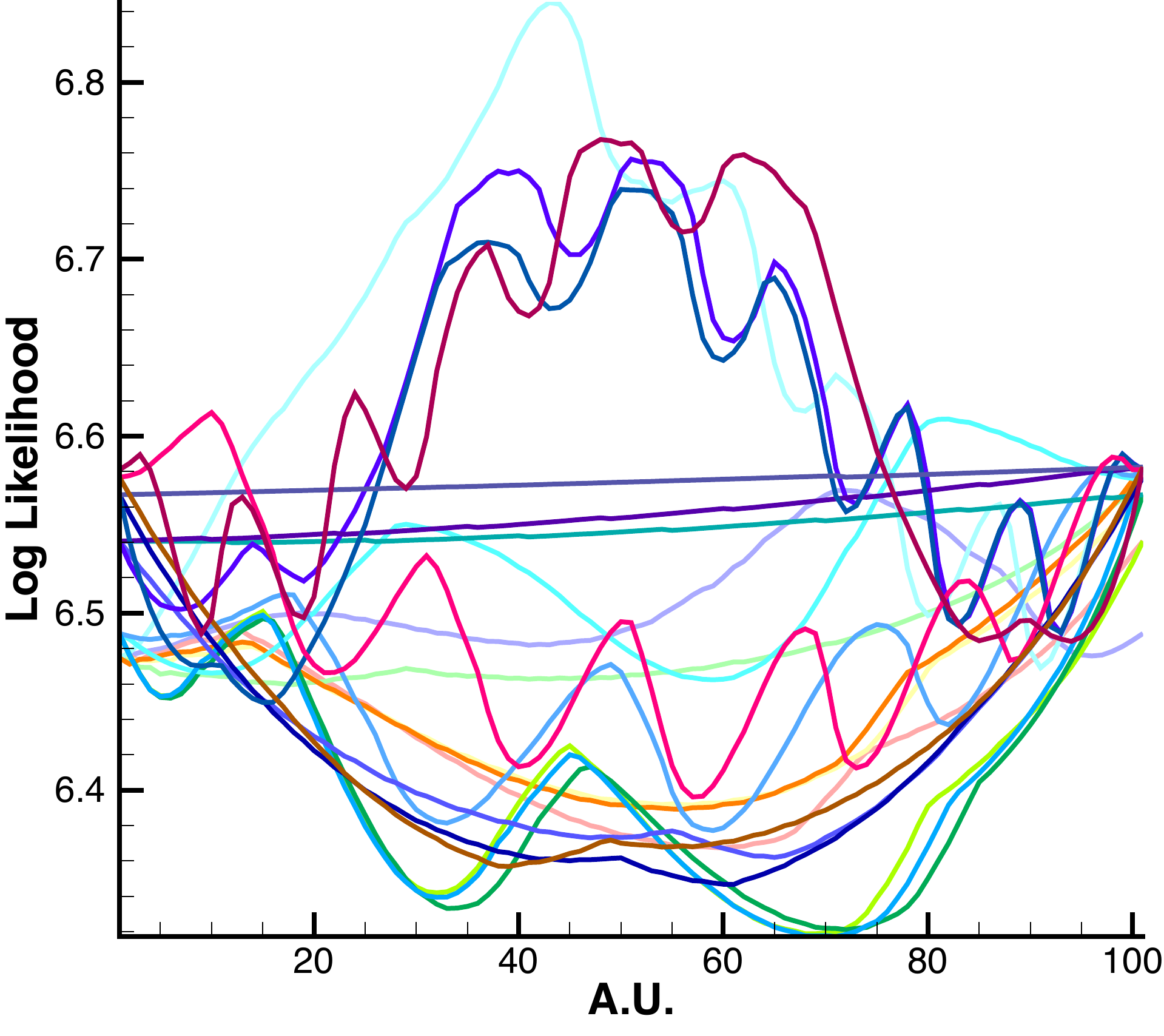}
\caption{
Negative logarithm of posterior distribution along
lines in parameter space.}
\label{fig:multimodal}
\end{figure}

Finally, we wish to point out that the triangle plot
suggests that the posterior has only one mode,
whereas our optimization found 12 ``distinct'' modes.
This apparent contradiction is caused by the modes
not being ``well separated''.
In fact we believe that the posterior distribution of this problem
is a nine dimensional analogue of the distribution illustrated in Figure~\ref{fig:IllustrateAlgorithm},
i.e., a plateau with several ``shallow dents''.
In this situation, the various modes may not be easy to spot in
a triangle plot because
(\textit{i}) we may not have sufficient resolution
due to limited sample size;
(\textit{ii}) the integration implied by viewing marginals only.
We can however, illustrate the complex structure of the posterior distribution
by evaluating the distribution along lines in parameter space
as in figure \ref{fig:multimodal}.
These lines are obtained as follows.
We take seven parameter vectors from our reference MCMC run
that result in the smallest negative logarithm of the posterior.
We then connect any two of these vectors by a line
and evaluate the posterior distribution at 100 equally spaced
points along these lines. 
We observe several ``dents'' along these lines.
These dents are deep enough to  
cause importance sampling algorithms that make
use of local geometry to fail.
For example, we tried implicit sampling (see, e.g., \cite{Chorin09,MorzfeldEtAl11}),
in which a proposal distribution is defined by the mode
and the Hessian of the negative log-posterior at the mode
and found that the method 
samples efficiently in the vicinity of the local minima,
but fails to explore the broader parameter space. 
The various local minima are also deep enough
to cause large integrated auto correlation times
in the MCMC runs (see above).

\section{Summary and conclusions}
\label{sec:Conclusions}
We presented iterative importance sampling schemes,
which use samples obtained by importance sampling
to construct a new proposal distribution that leads to more efficient importance sampling.
These ideas have been discussed before (see our literature review),
and the convergence of such iterations has also been studied.
It is important in practice that the iteration remains robust
when only few effective samples are produced at each step.
We studied this behavior in two test problems,
one in subsurface flow and one in combustion modeling.
In our implementation of ISA we made use of massively parallel computers 
and found that the iterations can robustly leverage parallelism,
which led to significant speed ups when compared to MCMC.
We also presented two strategies to initialize the iterations.
One strategy uses a ``short'' MCMC run,
the other makes use of optimization and Gaussian mixture models.

A natural question is: how do the computational requirements of ISA scale
with the dimension of the parameter vector?
Our numerical results provide little information here
because the dimension of the parameter vectors we considered are relatively small.
On the other hand, we argued in section~\ref{sec:ISAConvergence}
that the quality measure $R$ of Gaussian sampling problems
depends exponentially on the dimension of the parameter vector
(see also \cite{Owen,CLMMT16}).
This in turn implies that the number of samples required to achieve a given
number of ``effective samples'' also increases exponentially with dimension.
This may limit applicability of ISA, 
as the method would get stuck during the initialization phase
if $R$ were inevitably large in high dimensional problems.
On the other hand,
the Gaussian example that exhibits exponential dependence on
the dimension does not account for, e.g., sparsity of covariance matrices
of the target or proposal distributions.
If such sparsity were present, 
and if one could exploit sparsity during sampling,
then high dimensional problems may come within reach of ISA
or other importance sampling methods.
The situation is perhaps analogous to numerical linear algebra,
where computations with large dense matrices are infeasible,
while sparse matrices can be handled easily.
In fact, sparsity of covariances is exploited during ``covariance localization''
of ensemble Kalman filters and is the key to being able to effectively
solve estimation problems in numerical weather prediction
(where the state dimension is several hundred million).
If such strategies could be applied to parameter estimation problems,
e.g., by enforcing that some data only inform a subset of the 
parameters one wants to estimate,
then high-dimensional problems may become feasible.
In the applications presented here, 
we avoid this issue by selecting {\it {a priori}} only a small set of parameters.

\bibliographystyle{siamplain}
\bibliography{References}
\end{document}